\documentclass{article}

\usepackage{amssymb}
\usepackage{dsfont}

\newtheorem{lem}{Lemma}
\newtheorem{defi}[lem]{Definition}
\newtheorem{theo}[lem]{Theorem}
\newtheorem{prop}[lem]{Proposition}

\newtheorem{fait}[lem]{Fact}
\newtheorem{cor}[lem]{Corollary}

\newcommand{\proof}{\noindent{\bf Proof.}~}
\newcommand{\qed}{\ \hfill$\square$\bigskip}

\newcommand{\rk}{\hbox{\rm rk}\,}

\newcommand{\GL}{\hbox{\rm GL}\,}
\newcommand{\SL}{\hbox{\rm SL}\,}
\newcommand{\PSL}{\hbox{\rm PSL}\,}

\newcommand{\Aut}{\hbox{\rm Aut}\,}

\newcommand{\<}{\langle}

\renewcommand{\>}{\rangle}
\renewcommand{\o}{^\circ}

\newcommand{\C}{\mathds{C}}

\newcommand{\Z}{\mathds{Z}}

\title{Small groups of finite Morley rank with involutions}
\author{Adrien Deloro\footnote{Rutgers University} 
and Eric Jaligot\footnote{Corresponding author --- Universit\'e de Lyon, CNRS and 
Universit\'e Lyon 1}}

\date{June 21, 2008}

\begin{document}
\maketitle

\begin{abstract}
By analogy with Thompson's classification of nonsolvable finite $N$-groups, we 
classify groups of finite Morley rank with solvable local subgroups of even and of mixed 
type. We also consider miscellaneous aspects concerning ``small" groups of finite 
Morley rank of odd type. 
\end{abstract}

\section{Introduction} 

When they are present, involutions play a major role in the classification of 
infinite simple groups of finite Morley in a way much reminiscent of their use 
in the Classification of the Finite Simple Groups. 
This is at least the goal of Borovik's Program to tranfer 
as much as possible of arguments based on involutions 
from the finite case to the case of finite Morley rank. 
In both cases, most critical configurations occur when considering 
``small" groups, which are at the same time, and by far, the most difficult to handle. 
In the present paper we will consider the easiest cases among these remarkably 
difficult cases. 

Groups of finite Morley rank are equipped with a rudimentary notion of dimension 
on their first-order definable sets which behaves, or at least can be seen 
as, an abstract version of the Zariski dimension of algebraic varieties over algebraically 
closed fields. We refer to \cite{BorovikNesin(Book)94, AltBorCher(Book)} for the 
developments of the theory of groups of finite Morley rank and its links with 
finite group theory and algebraic group theory, which it encapsulates in a much more 
general and unified theory. Finite groups are exactly the groups of Morley rank $0$, 
and for algebraic groups over algebraically closed fields (with no additional 
structure) the Morley rank corresponds to the geometric Zariski dimension. 

It is known from \cite{BorovikBurdgesCherlin07} that a connected group 
of finite Morley has either trivial or infinite Sylow $2$-subgroups. Considering 
connected groups of finite Morley rank with involutions one has thus only to focus on 
groups with infinite 
Sylow $2$-subgroups. Groups without involutions lead to situations similar to that 
of the Feit-Thompson (Odd Order) Theorem in finite group theory, with 
no known infinite analog and actually different problems in this case 
\cite{Jaligot01}. 

The preliminary result for Borovik's Program in presence of involutions 
is the following. 

\begin{fait}\label{FactStruct2Syl}
{\bf \cite{BorovikPoizat90}}
Let $G$ be a group of finite Morley rank. Then Sylow $2$-subgroups of $G$ are conjugate 
and if $S$ is one of them, then $S\o$ is nilpotent and a central product, 
with finite intersection, of a $2$-torus and a $2$-unipotent subgroup. 
\end{fait}

As usual, {\em Sylow $p$-subgroups} are defined as the maximal $p$-subgroups. 
A {\em $2$-torus} is a divisible abelian $2$-group of finite Pr\"{u}fer rank, and a 
{\em $2$-unipotent} group is a definable connected (nilpotent) $2$-group of bounded 
exponent. 
These are similar in nature to the Sylow $2$-subgroup structure of an algebraic group 
over an algebraically closed field of characteristic different from and equal to 
$2$, respectively. Accordingly, we usually say that the group $G$ has the follwing type, 
depending on the nontriviality of the $2$-torus $T$ and/or of the $2$-unipotent 
subgroup $U$: 

\begin{center}
\begin{tabular}{c|cc}
                  & $U\neq 1$          & $U=1$  \\ \hline
$T\neq 1$   & Mixed                & Odd\\
$T=1$         & Even                 & Degenerate  \\
\end{tabular}
\end{center}

We note that the terminology ``Degenerate" has been adopted at a time when 
connected groups of finite Morley rank with finite but nontrivial Sylow $2$-subgroups 
were possibilities not entirely excluded, but since \cite{BorovikBurdgesCherlin07} it seems 
more informative to speak just of groups without involutions as far as connected groups 
are concerned. 

For some reasons related to a global unipotence theory in groups of finite Morley 
rank explained in \cite{JaligotSimple08, FreconJaligot07}, the following notation 
may be used. For a group $G$ of finite Morley rank, $\pi$ a set of primes and $p$ a prime, 
we let $U_{(\infty,0),\pi}(G)$ denote the subgroup of $G$ generated by definable decent 
tori of $G$ which coincide with the definable hull of their Hall $\pi$-subgroups, and 
$U_{(p,\infty)}(G)$ the subgroup of $G$ generated by its $p$-unipotent subgroups. 
In general, a {\em decent torus} is a definable divisible abelian group of finite Morley 
rank which coincides with the {\em definable hull}, that is the smallest definable 
containing subgroup \cite[\S I 2.3]{AltBorCher(Book)}, of its (divisible) torsion 
subgroup, and for an arbitrary prime $p$ a {\em $p$-unipotent group} 
is a definable connected {\em nilpotent} 
$p$-subgroup of bounded exponent. Subgroups of the form 
$U_{(\infty,0),\pi}(G)$ and $U_{(p,\infty)}(G)$, and more generally 
any subgroup generated by an arbitrary family of definable connected subgroups, 
are definable and connected by a well known application of Zilber's generation lemma 
\cite[\S5.4]{BorovikNesin(Book)94}. 

When considering only the prime $2$, this allows one to define naturally the 
{\em odd} and {\em even} parts of an arbitrary group of finite Morley rank. 
In particular, and as the letters ``T" and ``U" are reserved 
for ``torus" and ``unipotent" respectively, one uses as in \cite{AltBorCher(Book)} the 
following simpler notation for the odd and even parts of a group $G$ of finite Morley rank: 

$$T_{2}(G)=U_{(\infty,0),{\{2\}}}(G) \mbox{~and~} U_{2}(G)=U_{(2,\infty)}(G).$$ 

Any group $G$ of finite Morley rank has a {\em solvable radical}, that is a unique 
maximal definable solvable normal subgroup, usually denoted $R(G)$ 
\cite[Theorem 7.3]{BorovikNesin(Book)94}. Following the general 
classification of simple groups of finite Morley rank {\em of even type} as 
algebraic groups over algebraically closed fields 
\cite{AltBorCher(Book)}, the even part of an arbitrary group of finite 
Morley rank is best described as follows. 

\begin{fait}\label{CorGenClassificEvenType}
{\bf \cite[Proposition II 4.8 and Proposition X 1]{AltBorCher(Book)}}
Let $G$ be a group of finite Morley rank. Then $U_{2}(G)$, modulo 
its solvable radical, is a direct product of finitely many definable simple 
algebraic factors over algebraically closed fields of characteristic $2$. 
\end{fait}

The most general notion of ``smallness" for a group of finite Morley rank, incorporating 
notably all solvable groups and Chevalley groups of type $\PSL_{2}$ and $\SL_{2}$ over 
algebraically closed fields, is the following. 

\begin{defi}
A group of finite Morley rank is {\em locally$\o$ solvable$\o$} if the 
connected component of the normalizer of any infinite 
solvable subgroup is solvable. This is equivalent to require that the 
connected component of the normalizer of any nontrivial definable connected abelian 
is solvable by 
\cite[Lemma 3.4 $(4)$]{DeloroJaligotI}. 
\end{defi}

We refer 
to \cite{DeloroJaligotI} for a detailed study of such groups without any special assumption 
on the presence of involutions, and their analogies with those 
encountered in the Feit-Thompson Theorem and in Thompson's classification 
of finite ``$N$-groups" with involutions 
\cite{Thompson68, Thompson70, Thompson71, Thompson73}. 

We note that the latter series of papers also corresponds in the finite case to a transfer 
from {\em simplicity} to {\em nonsolvability} of certain arguments for 
``small" groups, and one of the goals of the present paper concerning groups of 
finite Morley rank is the same. 
In particular, we will study {\em nonsolvable} 
locally$\o$ solvable$\o$ groups of even and mixed types in 
Sections \ref{SectionEventype} and \ref{SectionMixedType} 
respectively, by using results and/or technics from the study of {\em simple} 
groups of finite Morley rank \cite{AltBorCher(Book)}. As naturally expectable, 
at least according to a long-standing feeling that groups of finite Morley rank resemble 
algebraic groups, our conclusion will be the following. 

\begin{theo}\label{TheoSmallGpsEvenAndMixedType}
Let $G$ be a locally$\o$ solvable$\o$ group of finite Morley rank with an infinite 
Sylow $2$-subgroup. Then exactly one of the following three cases occur. 
\begin{itemize}
\item[$(1)$]
$G\o$ is solvable. 
\item[$(2)$]
$G\o\simeq \PSL_{2}(K)$ for some algebraically closed field $K$ of characteristic $2$, 
in which case $G=G\o \times E$ fome some finite subgroup $E$. 
\item[$(3)$]
$G\o$ is nonsolvable and has odd type. 
\end{itemize}
\end{theo}

The case of nonsolvable connected groups of odd type has been studied 
in the thick series of consecutive works 
\cite{jal3prep, CherlinJaligot2004, BurdgesCherlinJaligot07, 
DeloroThesis, Deloro05, Deloro07}, more precisely in the simple case, implying in this 
process large portions of the current developments of the theory 
of groups of finite Morley rank. A kind of 
``final" version of this voluminous work will be found in 
\cite{DeloroJaligotII}. As a preparation, we will consider here certain 
specialized topics concerning groups of odd type, the case of {\em solvable} 
groups of odd type in Section \ref{SectionSolvGpsOddType} with generalities 
on involutive actions in Section \ref{SectionInvoActions}, 
and the case of groups of odd type with ``very small" Sylow 
$2$-subgroups in Section \ref{SectionGpsPruferRankOne}. 

In Section \ref{SectionBorovikCartan} we will also consider centralizers of involutions 
in groups of finite Morley rank when such involutions satisfy certain geometric properties 
reminiscent of small groups such as $\PSL_{2}$. This requires specific arguments 
analogous to the Cartan polar decomposition. 

\subsection{Locally$\o$ solvable$\o$ groups of even type}\label{SectionEventype}

As in the case of simple groups \cite{AltBorCher(Book)}, 
the elimination of connected nonsolvable groups of mixed type primarily depends on a 
classification of even type groups in the context of 
locally$\o$ solvable$\o$ groups. The most relevant statement is the following. 

\begin{theo}\label{MainTheoChar2}
Let $G$ be a locally$\o$ solvable$\o$ group of finite Morley rank of even type. 
Then exactly one of the following two cases occur. 
\begin{itemize}
\item[$(1)$]
$G\o$ is solvable, or 
\item[$(2)$]
$G\o \simeq \PSL_{2}(K)$ for some algebraically closed field $K$ of characteristic $2$, 
in which case $G=G\o \times E$ fome some finite subgroup $E$. 
\end{itemize} 
\end{theo}

\bigskip
We proceed to the proof of Theorem \ref{MainTheoChar2}. 
Let $G$ be a locally$\o$ solvable$\o$ group of finite Morley rank of even type, 
which may be assumed to be connected as long as one considers only its 
connected component. We then have 
$$R(U_{2}(G)) \trianglelefteq U_{2}(G) \trianglelefteq G.$$
If $R\o(U_{2}(G))$ is nontrivial, then $U_{2}(G)$ is solvable by local$\o$ solvability$\o$ 
of $G$. It follows then that $G\o$ is solvable by local$\o$ solvability$\o$ again, 
so we are in case $(1)$ of Theorem \ref{MainTheoChar2}. 

Assuming now that we are not in case $(1)$ of Theorem \ref{MainTheoChar2}, 
we have thus 
$R(U_{2}(G))$ finite. Dividing $U_{2}(G)$ by its finite 
solvable radical, one gets a semisimple group, which is still 
locally$\o$ solvable$\o$ by 
\cite[Lemma 3.5]{DeloroJaligotI}. 
Let $$H=U_{2}(G)/R(U_{2}(G)).$$ 
By Fact \ref{CorGenClassificEvenType}, $H$ is a direct 
product of finitely many definable normal simple subgroups. 
One sees that it is in fact a single definable normal simple subgroup by 
local$\o$ solvability$\o$ of $H$. Hence $H$ is a simple group of even type. 

By the classification of the simple groups of even type, the main 
theorem of \cite{AltBorCher(Book)}, 
$H$ is a simple algebraic group over an algebraically closed field of characteristic $2$. 
By local$\o$ solvability$\o$ of $H$, one concludes that 
$$H\simeq \PSL_{2}(K)\leqno(*)$$
for some algebraically closed field $K$ of characteristic $2$. The analysis can be 
continued as follows. 

\begin{lem}
$Z(U_{2}(G))=1$. 
\end{lem}
\proof
We have $Z(U_{2}(G))\leq R(U_{2}(G))$ which is finite, 
and in fact one has equality by 
\cite[Fact 3.14]{DeloroJaligotI}. 

Now $U_{2}(G)$ is a connected group. Its commutator subgroup is definable and connected 
by a well known corollary of Zilber's generation lemma 
\cite[Corollary 5.29]{BorovikNesin(Book)94}. As $U_{2}(G)/Z(U_{2}(G))=H$ 
is a nonabelian simple connected group, one gets that $U_{2}(G)'$ covers $H$, and thus 
$U_{2}(G)=Z(U_{2}(G))\cdot U_{2}(G)'$. Now finiteness of the center together 
with the connectedness of $U_{2}(G)$ forces that $U_{2}(G)$ is {\em perfect}, 
that is equal to its commutator subgroup. 

Then the result of \cite{AltinelCherlin99} on central extensions of algebraic groups 
implies that $U_{2}(G)$ is a Chevalley group over the same field $K$. 
As we are in characteristic $2$ and dealing with $\PSL_{2}$, we conclude that the 
center is trivial. 
\qed

\bigskip
At this point we have $U_{2}(G)\simeq \PSL_{2}(K)$ for some algebraically 
closed field $K$ of characteristic $2$. 

\begin{lem}
$G=U_{2}(G)$. 
\end{lem}
\proof
Let $U$ be a maximal $2$-unipotent subgroup of $U_{2}(G)$. 
By conjugacy of such subgroups in $U_{2}(G)\simeq \PSL_{2}$, a Frattini Argument gives 
$$G=U_{2}(G)\cdot N\o(U)$$
and the factor $N\o(U)$ is solvable by local$\o$ solvability$\o$. 
Now $N\o(U)$ acts on $U_{2}(G)\simeq \PSL_{2}(K)$ by inner automorphisms as 
$\PSL_{2}$ has no graph automorphisms. In particular 
$$N\o(U)=N_{U_{2}(G)}(U) \times C_{G}(U_{2}(G))$$
and as the latter factor is finite by local$\o$ solvability$\o$ of $G$ and nonsolvability 
of $U_{2}(G)\simeq \PSL_{2}$, one even has by connectedness of $N\o(U)$ that 
$N\o(U)=N_{U_{2}(G)}(U)\leq U_{2}(G)$, and it follows that 
$G=U_{2}(G)\cdot N\o(U)\leq U_{2}(G)$. 
\qed

\bigskip
We have thus shown that $G\simeq \PSL_{2}(K)$ for some algebraically closed field 
$K$ of characteristic $2$ whenever $G$ is connected and not solvable as in case $(1)$ of 
Theorem \ref{MainTheoChar2}. To complete the statement 
as in case $(2)$ of that theorem, it just remains to show the following. 

\begin{lem}
Let $G$ be a locally$\o$ solvable$\o$ group of finite Morley rank such that 
$G\o\simeq \PSL_{2}(K)$ for some algebraically closed field $K$ of characteristic $2$. 
Then $G=G\o \times E$ fome some finite subgroup $E$. 
\end{lem}
\proof
Again, by a Frattini Argument, $G=G\o \cdot N(U)$ for some maximal $2$-unipotent 
subgroup $U$ of $G\o$, and as there are no 
graph automorphisms of $\PSL_{2}$ we get 
$N(U)=N_{G\o}(U)\times C_{N(U)}(G\o)$, and thus 
$G=G\o\times C_{G}(G\o)$. Now $E=C_{G}(G\o)$ is the desired group. 
\qed

\bigskip
This completes the proof of Theorem \ref{MainTheoChar2}. 

The reader might however wonder whether one really needs the big gun of the full 
classification of simple groups of even type for the isomorphism $(*)$ above. 
Fortunately, one can obtain this isomorphism much more directly in the 
case of locally$\o$ solvable$\o$ groups. Here is how a baby version 
of the proof would work along the lines of the original papers in the simple case, 
even though we are not going to review all details entirely. 

We thus have a simple locally$\o$ solvable$\o$ group $H$ of finite Morley rank 
of even type, and we want to show that $H\simeq \PSL_{2}(K)$ for some 
algebraically closed field $K$ of characteristic $2$. 

Fix $U$ a maximal $2$-unipotent subgroup of $H$. 
Let $B=N\o(U)$ and $M=N(U)=N(B)$. We know that $B$ is a Borel subgroup of 
$H$ by 
\cite[Lemma 3.9]{DeloroJaligotI}. 
One sees easily that $M$ is {\em weakly embedded} in $H$, which means that it is a proper 
definable subgroup containing an infinite Sylow $2$-subgroup and such that 
$M\cap M^{h}$ has finite Sylow $2$-subgroups for any element $h$ of $H$ not in $M$ 
\cite{AltBorCher99}. Actually the strict inclusion $M<H$ follows from the 
nonsolvability of the ambient locally$\o$ solvable$\o$ group $H$ (see 
\cite[Lemma 3.7]{DeloroJaligotI}).
The property of finiteness of 
Sylow $2$-subgroups of intersections of distinct conjugates of $M$ follows from 
specific Uniqueness Theorems, analogous to those of the so-called 
{\em Bender Method} in finite group theory, available in the 
context of locally$\o$ solvable$\o$ groups of finite Morley rank, 
\cite[Corollary 4.4]{DeloroJaligotI} 
or 
\cite[Corollary 5.12]{DeloroJaligotI}. 

One sees easily with the same kind of arguments that for any nontrivial 
$2$-unipotent subgroup $V$ of $U$, $N\o(X)\leq M$ for any infinite 
definable subgroup $X$ of $C\o(V)$. This is because $N\o(X)$ is solvable 
by local$\o$ solvability$\o$, and thus one can use the Uniqueness Theorems of 
\cite[\S4.1]{DeloroJaligotI}, 
more specifically 
\cite[Corollary 4.4]{DeloroJaligotI}
according to which $B$ is the unique Borel subgroup containing any of its 
nontrivial $2$-unipotent subgroups. 

This shunts the most difficult part of the analysis, \cite[Th\'eor\`eme 4.1]{Jaligot01-a}, 
reducing essentially to the 
situation of \cite[\S3]{Jaligot01-a}, with $M\o=B$ solvable by local$\o$ solvability$\o$ 
(a rather undirect fact in the general case of \cite{AltBorCher(Book)}, but rather 
direct in the case of \cite{Jaligot01-a}). 
We leave to the reader the pleasure of accomplishing the final recognition of 
$\PSL_{2}$ along the lines of arguments and the computations of 
\cite[Th\'eor\`eme 3.1]{Jaligot01-a}, using here the fact that all normalizers$\o$ 
of nontrivial solvable infinite subgroups are solvable. 

\subsection{Locally$\o$ solvable$\o$ groups of mixed type}\label{SectionMixedType}

A corollary of the full classification of simple groups of finite Morley rank of even 
type \cite{AltBorCher(Book)} and of the arguments of \cite{Jaligot99} is the following. 

\begin{fait}
{\bf \cite{AltBorCher(Book)}}
There is no simple group of finite Morley rank of mixed type. 
\end{fait}

We obtain a similar result for connected locally$\o$ solvable$\o$ groups of finite 
Morley rank replacing the simplicity assumption by a mere nonsolvability 
assumption, which is best stated in the following form. 

\begin{theo}\label{MixedTypeTheorem}
Let $G$ be a locally$\o$ solvable$\o$ group of finite Morley rank of mixed type. 
Then $G\o$ is solvable. 
\end{theo}

Theorem \ref{MixedTypeTheorem} can be deduced as a special case of 
the general theory developed in \cite{AltBorCher(Book)}, and at the end of 
the present section we will review how this can be done. In any case, the core of the 
proof boils down in our special context to the argument below. 

To prove directly Theorem \ref{MixedTypeTheorem}, 
we assume toward a contradiction that $G$ is a connected counterexample of 
minimal rank to the above statement. 

As $G$ is not solvable, it has a finite solvable radical by local$\o$ solvability$\o$. 
Dividing by the latter, one gets a group still of mixed type, still 
locally$\o$ solvable$\o$ but now semisimple 
(\cite[Lemma 3.15]{DeloroJaligotI}), 
and in which all proper nonsolvable 
definable connected subgroups are not of mixed type by minimality. 

Fix $U$ a maximal $2$-unipotent subgroup of $G$. 
The first main claim is 
$$M:=N(U)\mbox{~is weakly embedded in~}G.$$
This is a special application of 
\cite[Corollary 5.12]{DeloroJaligotI} 
for the prime $p=2$. We also note that for this particular prime $p=2$ it is necessary that a 
nontrivial $2$-torus commutes with a nontrivial $2$-unipotent subgroup by 
Fact \ref{FactStruct2Syl}, so that the general strategy developed in the 
simple case \cite{AltBorCher97} also works here, as explained after 
\cite[Corollary 5.12]{DeloroJaligotI}. 
It provides in particular a proof very similar to 
the one used in the simple case: $M=N(U)=N(U^{\perp})$ where 
$U^{\perp}=T_{2}(C(U))$, and one checks easily that $M$ contains the normalizer of 
each of its nontrivial $2$-unipotent subgroups (By the Uniqueness Theorem of 
\cite{DeloroJaligotI}) and similarly for its $2$-tori $T$: $N\o(T)$ is solvable by 
local$\o$ solvability$\o$, contains $U$, and the Uniqueness Theorem applies 
again. 

The next point is the following remark. 

\begin{fait}\label{CorNSControlFusion}
{\bf \cite[Fait 2.18]{Jaligot99}}
Let $G$ be a group of finite Morley rank, $S$ a Sylow $2$-subgroup of $G$, $T$ 
the maximal $2$-torus of $S\o$, and $t$ an element of $T$. Then 
$t^{G}\cap S\o$ is contained in $T$, and is in particular finite. 
\end{fait}
\proof
The first claim follows from an argument of control of fusion in $p$-tori by their 
normalizers, which has been known for a long time for the particular prime $p=2$ 
\cite[Lemma 10.22]{BorovikNesin(Book)94}, 
and the present formulation can be tracked in \cite[Fact 2.48]{Altinel96}. 
Anyway we refer to 
\cite[Corollary 2.20]{DeloroJaligotI} 
for the most general formulation of such arguments of control of fusion, 
in a form which directly applies here. 

The finiteness of $t^{G}\cap S\o$ follows, as the $2$-torus $T$ has only 
finitely many elements of order $2^{n}$ for each $n$ \cite{BorovikPoizat90}. 
\qed

\bigskip
A proper definable subgroup $M$ is {\em strongly embedded} if it has nontrivial 
Sylow $2$-subgroups and $M\cap M^{g}$ has trivial Sylow $2$-subgroups for 
any element $g$ of $G$ not in $M$. There is a much 
similar notion in the finite case, used notably by Bender, and the notion of weak 
embedding is its neoclassical revival in the case of groups of finite Morley rank. 

The next point is then the following. 

\begin{lem}
$M$ is not strongly embedded in $G$. 
\end{lem}
\proof
Assume towards a contradiction $M$ strongly embedded. 
Then its involutions would necessarily be conjugate 
\cite[Theorem 10.19]{BorovikNesin(Book)94}. 
In particular an involution of a maximal $2$-torus $T$ of $M$ would be conjugate 
to any involution of $U$. But there are infinitely many involutions in $2$-unipotent 
groups \cite{BorovikPoizat90}, and this 
contradicts Fact \ref{CorNSControlFusion}. 
\qed

\bigskip
By general results on strong/weak embedding (see \cite[Part B]{AltBorCher(Book)}), 
one then concludes that there exists an involution $\alpha$ of $M$ which 
is {\em problematic} in the sense that 
$$C(\alpha) \nleq M$$
and the next step consists in applying Theorem \ref{MainTheoChar2} to its 
centralizer$\o$. 

\begin{lem}\label{LemalphaProbCalphaPSL2}
Let $\alpha$ be a problematic involution of $M$. Then 
$C\o(\alpha)\simeq \PSL_{2}(K)$ for some algebraically closed field $K$ of 
characteristic $2$. 
\end{lem}
\proof
The involution $\alpha$ normalizes $U=U_{2}(B)$, as connected solvable 
groups of finite Morley rank have only one maximal $p$-unipotent subgroup 
(see for example 
\cite[Fact 2.15]{DeloroJaligotI} 
for the most general and contemporary version of this, 
or \cite{Nesin90-b} for the oldest results from which it can be deduced). 
As infinite nilpotent-by-finite $p$-groups of finite Morley rank contain 
infinitely many central elements of order $p$ 
\cite{BorovikPoizat90}, $\alpha$ centralizes a nontrivial $2$-unipotent subgroup of $U$. 

We claim that $C\o(\alpha)$ is not solvable. This follows from the 
Uniqueness Theorems of 
\cite[\S4.1]{DeloroJaligotI}, 
according to which a 
nontrivial $p$-unipotent subgroup of any locally$\o$ solvable$\o$ group of finite Morley 
rank is contained in a unique Borel subgroup. In particular $B$ is the unique 
Borel subgroup containing $C\o_{U}(\alpha)$. 
Assuming $C\o(\alpha)$ solvable, we would get in particular 
$C\o(\alpha)\leq B$, and $C(\alpha)\leq N(C\o_{U}(\alpha))\leq N(B)=M$, 
as $B$ is the unique Borel subgroup containing $C\o_{U}(\alpha)$, a contradiction since 
$\alpha$ is problematic. 

By semisimplicity of $G$, $C\o(\alpha)<G$, 
and by minimality, $C\o(\alpha)$ cannot be of mixed type. 
As it contains a nontrivial $2$-unipotent subgroup as just seen, $C\o(\alpha)$ is a 
locally$\o$ solvable$\o$ group of even type. Now 
Theorem \ref{MainTheoChar2} yields the desired isomorphism type of 
$C\o(\alpha)$. 
\qed

\bigskip
For the final step we are now in a position to conclude as at the end of 
\cite{Jaligot99}, using the same relevant technical lemma on involutions 
isolated to tackle the configuration appearing. 

\begin{fait}\label{FactJaligot99}
{\bf \cite[Lemme 4.1]{Jaligot99}}
Let $G$ be a group of finite Morley rank, with involutions 
$i$, $t$, $\alpha$, and $\alpha'$ satisfying the following five conditions. 
\begin{itemize}
\item[$(1)$]
$i$ and $t$ are not conjugate. 
\item[$(2)$]
$U_{2}(C(\alpha))\simeq \PSL_{2}(K)$ for some algebraically closed 
field $K$ of characteristic $2$. 
\item[$(3)$]
$\alpha'$ is the unique involution of the definable hull $H(it)$ of $it$ 
($i\cdot t$, not it). 
\item[$(4)$]
$i\in U_{2}(C(\alpha))$, $t\in C(\alpha)$. 
\item[$(5)$]
$\alpha' \notin U_{2}(C(\alpha))$. 
\end{itemize}
Then $t\alpha' \in U_{2}(C(\alpha))$. 
\end{fait}

\bigskip
\noindent
{\bf Proof of Theorem \ref{MixedTypeTheorem}.}
We keep the previous notations. 

Applying 
\cite[Lemma 3.35 $(2)$]{DeloroJaligotI}
to the conjugacy 
class of a $2$-toral involution $t$ of $G$, one gets as $G$ is not solvable that 
$t^{G}\cap M$ is not generic in $t^{G}$ (notice that $t^{G}$ is infinite as 
$C(t)<G$ and $G$ is connected). In particular one may assume, choosing a suitable 
conjugate of $t$, that $t$ is a $2$-toral involution not in $M$. 

Let $\alpha_{0}$ be a problematic involution of $M$. As seen in the proof of 
Lemma \ref{LemalphaProbCalphaPSL2}, there 
exists an involution $i_{0}$ in $C\o_{U}(\alpha_{0})$. As involutions of the latter 
group are all conjugate in $C\o(\alpha_{0})$ by Lemma \ref{LemalphaProbCalphaPSL2}, 
Corollary \ref{CorNSControlFusion} implies that $t$ and $i_{0}$ are not conjugate. 
In particular there exists an involution $\alpha$ in the definable hull 
$H(i_{0}t)$ of $i_{0}t$ \cite[Proposition 10.2]{BorovikNesin(Book)94}. 
We have 
$$[i_{0},\alpha]=1 \mbox{~and~} [t,\alpha]=1.$$ 
Lemma \ref{LemalphaProbCalphaPSL2} implies in particular that 
problematic involutions of $M$ never belong to a unipotent subgroup, as they cannot 
centralize a nontrivial $2$-torus. In particular $i_{0}$ is not problematic in $M$, 
i.e., $C(i_{0})\leq M$. The first commutation implies thus that 
$\alpha$ is in $M$, and as $t\notin M$ the second commutation implies that 
$\alpha$ is a problematic involution of $M$. 

Take an involution $i$ in $C\o_{U}(\alpha)$. Again $i$ and $t$ are not conjugate, 
hence there is an involution $\alpha'$ in $H(it)$. Similarly, $\alpha'$ commutes 
with $i$ and $t$, and is thus a problematic involution of $M$. 

One can check that the five conditions of Fact \ref{FactJaligot99} are met for these 
involutions $i$, $t$, $\alpha$, and $\alpha'$. For the third for example, 
we note that $H(it)$ contains no nontrivial $2$-torus, as $\alpha' \in H(it)$ and 
by the structure of centralizers of problematic involutions. For the last point, 
this is by the fact that problematic involutions cannot belong to a $2$-unipotent 
subgroup. The conclusion of Fact \ref{FactJaligot99} is then 
$$t\alpha' \in U_{2}(C(\alpha)).$$
As $t$ centralizes $\alpha$, it acts on $C\o(\alpha)$, and by the structure of the latter 
the action must be by inner automorphism. Since $t$ centralizes $t\alpha'$, it centralizes 
the Sylow $2$-subgroup $A_{t\alpha'}$ of $C\o(\alpha)$ containing $t\alpha'$. 

Now $A_{t\alpha'}\leq C(t)<G$, and hence by our minimality assumption 
$C(t)$ is of even type. But $t$ belongs to a nontrivial $2$-torus, and 
in particular centralizes it. This is a contradiction which ends the proof of 
Theorem \ref{MixedTypeTheorem}. 
\qed

\bigskip
Another approach to Theorem \ref{MixedTypeTheorem} would be to use the fact that 
$$[T_{2}(G),U_{2}(G)]=1.$$ 
This is a general corollary of the full 
classification of simple groups of finite Morley rank of even type 
\cite[Lemma V 2.3, Theorem X 1, Chap. V Mixed Type $L^{*}$ Theorem]{AltBorCher(Book)}. 
With this commutation, one concludes easily to the solvability of $G$ when 
both factors in the commutator are nontrivial, as both the normal factors are 
then solvable by local$\o$ solvability$\o$ of $G$. If one tries to prove directly 
this commutation in our case, then as in 
\cite[Lemma V 2.3]{AltBorCher(Book)} one may want to look at the action of 
$2$-tori on $U_{2}(G)$. As in Section \ref{SectionEventype}, 
or \cite[Proposition II 6.2, Lemma II 6.3]{AltBorCher(Book)}, one is then interested 
into the socle of $U_{2}(G)$ modulo $R(U_{2}(G))$, a direct product of 
connected simple factors. If no such factor is of mixed type, then one can conclude 
that $[T_{2}(G),U_{2}(G)]=1$ as in \cite[Lemma V 2.3]{AltBorCher(Book)}. 
Otherwise, the situation reduces to the case where $G$ is a simple group such that 
$G=U_{2}(G)=T_{2}(G)$, and in our case this was disposed of by the core of the proof 
given above. 

As it is was explained in 
\cite[\S3.3]{DeloroJaligotI}, 
there might be serious obstructions 
if one wants a version of Theorem \ref{MixedTypeTheorem} for primes different from $2$. 

We now note that Theorem \ref{TheoSmallGpsEvenAndMixedType} 
follows from Theorems \ref{MainTheoChar2} and \ref{MixedTypeTheorem}. 
By-the-by, we mention the following corollary of 
Theorem \ref{TheoSmallGpsEvenAndMixedType} on connectedness 
of Sylow $2$-subgroups in small groups of finite Morley rank. 

\begin{cor}
Let $G$ be a connected locally$\o$ solvable$\o$ group of finite Morley rank 
with a nonconnected Sylow $2$-subgroup. Then $G$ is nonsolvable and of odd type. 
\end{cor}
\proof
It is known that Hall $\pi$-subgroups are connected in any connected solvable 
group of finite Morley rank and for any set $\pi$ of primes. In $\PSL_{2}$ over 
some algebraically closed field $K$ of characteristic $2$, Sylow $2$-subgroups are 
also connected. It can be seen either by the transitivity of the action on them of 
their normalizers or by the fact that they are definably isomorphic to the additive 
group of the ground field $K$, which is interpretable and of finite Morley rank. 

It follows that the only remaining possibility in 
Theorem \ref{TheoSmallGpsEvenAndMixedType} 
is case $(3)$ in that theorem. 
\qed

\subsection{Involutive actions}\label{SectionInvoActions}

The structure of nonsolvable locally$\o$ solvable$\o$ groups of odd type, 
left undetermined in Theorem \ref{TheoSmallGpsEvenAndMixedType}, 
will be considered in \cite{DeloroJaligotII} as mentioned in the 
introduction. Now we merely deal with specialized topics concerning groups 
of odd type in general and which will in particular be applied in the lengthy 
analysis of \cite{DeloroJaligotII}. 

Before moving ahead, we concentrate in general on definable involutive automorphisms of 
groups of finite Morley rank. When $\alpha$ is an automorphism of a group $G$, we let 
$$G^{+_{\alpha}}=\{g\in G~|~g^\alpha=g\}
\mbox{~and~} 
G^{-_{\alpha}}=\{g\in G~|~g^\alpha=g^{-1}\}.$$
Both sets are definable whenever $\alpha$ is a definable automorphism of $G$. 
In general, only the centralizer $G^{+_{\alpha}}$ of $\alpha$ in $G$ is necessarily a 
subgroup of $G$. When there is no risk of confusion 
between different possible automorphisms $\alpha$, we will sometimes omit the 
subscript $\cdot_{\alpha}$ in the notation of the two sets as above, and thus just speak 
of $G^{+}$ and $G^{-}$. 

We start by mentioning the older results on involutive definable automorphisms 
$\alpha$ of groups of finite Morley rank in terms of $G^{+}$ and $G^{-}$. 

\begin{fait}\label{FaitActionInvCentFiniG0inv}
{\bf \cite{Nesin90-a}}
Let $\alpha$ be a definable involutive automorphism of a group $G$ of finite 
Morley rank. If $\alpha$ fixes only finitely many elements in $G$, then $G$ has 
a definable (abelian) subgroup of finite index inverted by $\alpha$, i.e., 
$G\o\subseteq G^{-}$. 
\end{fait}

\begin{fait}\label{ActionInvOn2PerpGroup}
{\bf \cite[Ex. 14, p. 73]{BorovikNesin(Book)94}}
Let $G$ be a group of finite Morley rank without involutions and $\alpha$ a 
definable involutive automorphism of $G$. Then $G^{-}$ is $2$-divisible and 
one has a decomposition $G=G^{+}\cdot G^{-}$, where 
the corresponding multiplication map is one-to-one. 
\end{fait}

In the present section we are essentially going to give a generalization of 
Fact \ref{ActionInvOn2PerpGroup} when the group $G$ contains a central and 
divisible Sylow $2$-subgroup. We first note the following general lemma. 

\begin{lem}\label{LemGCentpToreGpDiv}
Let $G$ be a group of finite Morley rank whose Sylow $p$-subgroup is a central $p$-torus. 
Then $G$ is $p$-divisible.
\end{lem}
\proof
Let $g$ be an arbitrary element of $G$. Then the definable hull $H(g)$ of the cyclic 
group $\<g\>$ can be written as 
$$H(g) = \Delta \oplus \< \zeta \>$$ 
where $\Delta$ is $p$-divisible and $\zeta$ a $p$-element, by using 
\cite[Lemma 2.16]{AltBorCher(Book)} and by decomposing the finite cyclic part into 
its $p$-primary component and a (cyclic) complement. 
Let $h$ be a $p$-th root of 
$g\zeta^{-1}$ in $\Delta$. Let $\eta$ be a $p$-th root of $\zeta$ in the 
Sylow $p$-subgroup. As $\eta$ is central in $G$, one gets 
$$g=g\zeta^{-1}\zeta=h^p\eta^p=(h\eta)^p.$$
\qed

In general, a group of even type is not $2$-divisible. One may wonder whether 
groups of odd type are $2$-divisible, but this need not hold in general neither as the following 
example shows. If $B$ is a Borel subgroup of $\SL_2(K)$, with $K$ an algebraically closed 
field of characteristic different from $2$, $u$ a nontrivial unipotent element in $B$, and 
$i$ the unique (central) involution of $\SL_{2}(K)$, then $ui$ is not a square in 
$\SL_{2}(K)$. In particular $B$ is a connected $2$-step solvable group of finite 
Morley rank of odd type, and $B$ is not $2$-divisible! 

Our generalization of Fact \ref{ActionInvOn2PerpGroup} is the following. 

\begin{theo}\label{actionalpha}
Let $G$ be a group of finite Morley rank whose Sylow $2$-subgroup is a 
(possibly trivial) central $2$-torus, and $\alpha$ a definable involutive 
automorphism of $G$. Then $G=G^+\cdot G^-$ where the fibers of the associated 
product map are finite. 
Furthermore one also has $G=(G^+)\o \cdot G^-$ whenever $G$ is connected. 
\end{theo}

The proof of Theorem \ref{actionalpha} generalizes that 
of Fact \ref{ActionInvOn2PerpGroup}. 
We note that Theorem \ref{actionalpha} also incorporates 
\cite[Fait 1.19]{Deloro07}. 
We shall need the following intermediate general lemma whose proof has a flavor 
similar to that of \cite[Lemme 4.45]{Jaligot01-a}. 

\begin{lem}\label{Cpm}
Let $G$ be a group of finite Morley rank whose Sylow $2$-subgroup is central. 
If $a^b = a^{-1}$ for two elements $a$ and $b$ of $G$, then $a$ is either the identity 
or an involution of $G$. 
\end{lem}
\proof
The set $C^{\pm}(a)=\{ x\in G, a^x = a^{\pm 1}\}$ is a definable subgroup of $G$. 
As $C^{\pm}(a)/C(a)$ has exponent at most $2$, there is a (possibly trivial) $2$-element in 
any coset of $C(a)$ in $C^{\pm}(a)$ by \cite[Lemma 2.18]{AltBorCher(Book)}. 
As $2$-elements are central by assumption, 
this proves that any coset of $C(a)$ in $C^{\pm}(a)$ is in $C(a)$, and thus $C^{\pm}(a)=C(a)$. 
In particular $a=a^{-1}$, and $a^2=1$.
\qed

\bigskip
\noindent
{\bf Proof of Theorem \ref{actionalpha}.} 
Recall that $G$ is a group of finite Morley rank whose Sylow $2$-subgroup is a 
(possibly trivial) central $2$-torus, and that $\alpha$ is an automorphism of $G$ 
of order at most $2$.

\bigskip
\noindent
{\bf Step 1.} 
{\it  
If $[g,\alpha]$ is a non-trivial $2$-element, then there exists $h$ in $G$ such that 
$[g,\alpha]=[h,\alpha]^2$.}

\medskip
\noindent
{\bf Proof.}
Assume $\zeta = [g, \alpha]$ has order $n=2^k$ for some $k$. 
As $G$ is $2$-divisible by Lemma \ref{LemGCentpToreGpDiv}, there exists an element 
$h$ in $G$ such that $h^{2}=g$. Now 
$\zeta = [g, \alpha] = [h^2, \alpha] = [h, \alpha]^h [h, \alpha]$ is central and has 
order $n$, so $h$ inverts $q = [h, \alpha]^n$. By Lemma \ref{Cpm}, $q$ has order 
at most $2$. In particular $[h, \alpha]$ is a $2$-element, and it is in particular central 
in $G$. One then gets 
$[g,\alpha]=[h,\alpha]^h [h,\alpha]=[h,\alpha]^2$. 
\qed

The moral of Step 1 is that in general $G^{-}$ needs not be $2$-divisible as in 
Fact \ref{ActionInvOn2PerpGroup}, but its subset 
$\{[g,\alpha], g\in G\}$ is $2$-divisible as we will see in the course of the next step. 

\bigskip
\noindent
{\bf Step 2.}
{\it $G = G^+ \cdot G^-$.}

\medskip
\noindent
{\bf Proof.}
Let $g \in G$. One can write the definable hull $H([g,\alpha])$ of $[g,\alpha]$ as 
in Lemma \ref{LemGCentpToreGpDiv} as 
$H([g,\alpha])=\Delta \oplus \< \zeta\>$, where 
$\Delta$ is a $2$-divisible group and $\zeta$ is a $2$-element. 
Let $h$ be a square root of $[g,\alpha]\zeta^{-1}$ in $\Delta$. 
Notice that $h$ centralizes $[g,\alpha]$ and is inverted by $\alpha$. Hence 
$[gh,\alpha]=[g,\alpha]^{h}[h,\alpha]=[g,\alpha] h^{-2}=\zeta$. 
By Step 1, there is $x \in G$ such that $\zeta = [x, \alpha]^2$. So $\eta = [x, \alpha]$ 
is a $2$-element inverted by $\alpha$, whose square is $\zeta$. 
As $\eta \in Z(G)$, we find that $\eta^{-1}h^{-1}$ is inverted by $\alpha$. Moreover, 
$[gh\eta, \alpha] = [g, \alpha]^{h\eta}[h, \alpha]^\eta[\eta, \alpha] = 
h^2 \zeta h^{-2} \eta^{-2} = 1$, so $g = (gh\eta)(\eta^{-1}h^{-1})$ 
is a suitable decomposition. 
\qed

\noindent
{\bf Step 3.}
{\it If $ax = by$ with $a$, $b\in G^{+}$ and $x$, $y\in G^{-}$, 
then $a^{-1}b$ is the identity or an involution fixed by $\alpha$. 
Hence the fibers of the product decomposition as in Step 2 are of 
cardinal at most $|\{k\in G~|~k^{2}=1\}|$.}

\medskip
\noindent
{\bf Proof.} 
We have 
$$(a^{-1}b)^y = (xy^{-1})^y = y^{-1} x = \alpha(yx^{-1}) =
\alpha(b^{-1} a) = b^{-1} a = (a^{-1}b)^{-1},$$ so $y$ inverts $a^{-1}b$.  
By Lemma \ref{Cpm}, $a^{-1}b$ has order at most $2$. 
In particular, $a^{-1}b$ lies in the central elementary abelian $2$-subgroup of $G$, 
and it is fixed by $\alpha$. We also note that the central elementary abelian $2$-subgroup of 
$G$ is exactly $\{k\in G~|~k^{2}=1\}$. 
\qed

\noindent
{\bf Step 4.}
{\it Left $G^+$-translates of the set $(G^+)\o \cdot G^-$ are disjoint or equal.}

\medskip
\noindent
{\bf Proof.} 
Assume that $a(G^+)\o \cdot G^-$ meets $b(G^+)\o \cdot G^-$, in say 
$ag_{+}g_{-}= bh_{+}h_{-}$ with natural notations. 
By Step 3, $z=(ag_{+})^{-1}(bh_{+})$ is central and fixed by $\alpha$. So 
$z=(bh_+)(ag_+)^{-1}$. Hence for any $b\gamma_+ \gamma_- \in b(G^+)\o \cdot G^-$ 
one finds 
$b\gamma_{+}\gamma_{-}=zz^{-1}b\gamma_{+}\gamma_{-}=a(g_+ h_+^{-1}\gamma_+) (\gamma_- z)$, 
which lies in $a (G^+)\o \cdot G^-$. The converse inclusion holds too. 
\qed

\noindent
{\bf Step 5.}
{\it Exactly $\deg(G)$ left $G^+$-translates of $(G^+)\o\cdot G^-$ cover $G$. 
In particular, if $G$ is connected, then $G=(G^+)\o G^-$.}

\medskip
\noindent
{\bf Proof.} 
We consider such left translates. They all have rank $\rk(G)$ by Step 3. 
As they are disjoint or equal by Step 4, exactly $\deg(G)$ of them suffice to cover $G$.
\qed

This ends the proof of Theorem \ref{actionalpha}.
\qed

\bigskip
The following results are not used here, but it is worth mentioning them for the 
sake of completeness of the present section. The first one is an important 
commutation principle, and the second one is a downward invariance result. 

\begin{fait}\label{HKK2divActionInv}
{\bf \cite[Lemme 3.1]{Deloro05}} 
Let $G$ be a group, $H$ and $K\leq N(H)$ two subgroups with $K$ $2$-divisible. 
Suppose that there is an involution $i$ in $G$ which inverts $K$ and centralizes or inverts 
$H$. Then $[H,K]=1$. 
\end{fait}

\begin{fait}\label{Snakes1}
{\bf (Compare with \cite[Fait 3.12]{Deloro05})}
Let $G$ be a group of finite Morley rank, $K\leq G$ a definable $2$-subgroup 
normalizing a definable subgroup $H$ of $G$, and $S$ a Sylow $2$-subgroup of $H$. 
Then an $H$-conjugate of $K$ normalizes $S$. 
\end{fait}

\proof
Consider the definable group $HK$, and let $\hat{S}$ denote a Sylow $2$-subgroup of 
$HK$ containing $S$. By conjugacy of Sylow $2$-subgroups, $K^{\gamma}\leq \hat{S}$ 
for some $\gamma$ in $HK$, and in fact $\gamma$ may be chosen in $H$. Now the 
$H$-conjugate $K^{\gamma}$ of $K$ normalizes $\hat{S}\cap H=S$.  
\qed

\subsection{Solvable groups of odd type}\label{SectionSolvGpsOddType}

Since the original unpublished \cite{jal3prep} on small groups 
of odd type, most of the work has consisted in developing certain arguments 
based on considerations involving {\em strongly real} elements, that is products of two 
involutions. The goal is ultimately to adapt to odd type groups such arguments 
from finite group theory, first imported by Nesin to groups of finite Morley rank. 

In $\PSL_{2}(K)$, with $K$ an algebraically closed field of characteristic different from 
$2$, the standard Borel subgroup contains two kinds of strongly real elements. Those 
of maximal unipotent subgroups on the one hand are inverted by 
involutions inside their respective Borel subgroups, and those in maximal tori on the other hand 
are inverted by involutions outside the Borel subgroup, and corresponding to 
liftings of elements of the Weyl group. We may call these 
strongly real elements ``insiders" and ``outsiders" respectively. 

When working with small groups of odd type, most complications arise when trying 
to get a control on outsiders, that is strongly real 
elements inside a Borel subgroup but such that the two involutions forming the product 
are outside this Borel subgroup. This is typically done for a {\em standard} Borel 
subgroup, that is a Borel subgroup $B$ containing the centralizer$\o$ of a particular 
involution. In this case outsiders of this Borel subgroup $B$, inverted by a fixed involution, 
are expected to form a torus of $B$ in the expected group $\PSL_{2}$. 
Other complications arise when one has to compare insiders and outsiders 
of such a Borel subgroup, and in this case one also occasionally needs in this delicate work 
a good understanding of insiders. We refer to \cite{DeloroJaligotII} anyway. 

The study of insiders merely boils down to 
the study of connected solvable groups of odd type, and this is the purpose of the present 
section. 
Recall first from \cite{BorovikPoizat90} that the connected 
component of the Sylow $2$-subgroup of a group of finite Morley rank is 
a direct product of {\em finitely many} copies of the quasicyclic group 
$\Z_{2^{\infty}}$, with this finite number called the 
{\em Pr\"{u}fer $2$-rank} of the group. As connected solvable groups of finite 
Morley rank have connected Sylow $2$-subgroups, the Sylow $2$-subgroup $S$ of 
a connected solvable group of odd type satisfies 
$$S=S\o \simeq {\Z_{2^{\infty}} \times \cdots \times \Z_{2^{\infty}}}$$
where the number of factors involved is the Pr\"{u}fer $2$-rank. 

Before moving ahead, we recall some more general background. 
As for the solvable radical, there is a similar notion of {\em Fitting} subgroup in any group 
$G$ of finite Morley rank \cite[Theorem 7.3]{BorovikNesin(Book)94}, where nilpotence 
replaces solvability, usually denoted by $F(G)$. A {\em Carter} subgroup of a group 
of finite Morley rank is a definable connected nilpotent subgroup of finite 
index in its normalizer, and it is nontrivial but a fact that any group of finite 
Morley rank contains a Carter subgroup \cite{JaligotFrecon}. 

\begin{lem}\label{LempToralCarter}
Let $G$ be a connected group of finite Morley rank without nontrivial $p$-unipotent 
subgroups for some prime $p$, and $t$ a $p$-element of $G$. Then 
\begin{itemize}
\item[$(1)$]
The element $t$ belongs to a Carter subgroup $Q$ of $G$, and $Q\leq C\o(t)$. 
\item[$(2)$]
If furthermore $G=NQ$ for some normal definable subgroup $N$, then 
$G=C\o(t)N$. 
\end{itemize}
\end{lem}
\proof
$(1)$. By the main result of \cite{BurdgesCherlinSemisimpleTorsion}, $t$ 
belongs to a $p$-torus, say $T_{p}$. Now any decent torus is contained in a Carter 
subgroup by the construction given in \cite{JaligotFrecon}. We have thus 
$t\in T_{p}\leq Q$ for some Carter subgroup $Q$ of $G$, and as any 
connected nilpotent group of finite Morley rank has a unique decent torus which is 
central, by earlier work of Nesin, one concludes that 
$Q\leq C\o(T_{p})\leq C\o(t)$. The reader can also consult 
\cite[Fact 4]{JaligotGenerixCosets} for similar facts in this direction. 

$(2)$. We have $G=QN\leq C\o(t)N\leq G$, and thus $G=C\o(t)N$. 
\qed 

\bigskip
Here are the main lemmas used in \cite{DeloroJaligotII} concerning 
strongly real elements inside connected solvable groups of odd type. 
We stress the fact that the Pr\"{u}fer $2$-rank is not $1$ in general, a technical 
complication which has to be entirely taken into account in 
\cite{Deloro07} and in \cite{DeloroJaligotII}. 
For an element normalizing a subgroup we confound below the element with 
the automorphism it induces by conjugation in the notation introduced at the begining of 
Section \ref{SectionInvoActions}. 

\begin{lem}[Inner Computation]\label{computationinside}
Let $B$ be a connected solvable group of finite Morley rank of odd type and 
let $j$ be an involution of $B$. Then 
\begin{itemize}
\item[$(1)$]
$B=C\o_B(j) \cdot [F\o(B)]^{-_j}$, where the fibers of the associated product 
map are finite. In particular $\rk(B)=\rk(C\o_B(j))+\rk([F\o(B)]^{-_j})$. 
\item[$(2)$]
For any definable subgroup $U$ of $Z(F\o(B))$ normal in $B$, 
the subgroup $U^{-_{j}}$ is normal in $B$. 
\end{itemize}
\end{lem}

\proof 
$(1)$. 
As $B$ is of odd type, there are finitely many involutions in the connected 
nilpotent group $F\o (B)$ \cite[\S6.2]{BorovikNesin(Book)94}. 

We claim that the fibers of the multiplication map 
$$P~:~C\o(j) \times (F\o(B))^{-_j} \rightarrow B$$ 
are finite. 
If $cf=c'f'$, with obvious notations, then $c'^{-1}c=f'f^{-1}$, and this element 
is centralized by $j$. Hence one gets $f'f^{-1}=(f'f^{-1})^{j}=f'^{-1}f$, and 
$f'^{2}=f^{2}$. The nilpotent group $F\o(B)$ has the form $T\ast U$, a central 
product with finite intersection, where $T$ denotes the maximal decent torus of 
$F\o(B)$ and $U$ is a definable connected nilpotent subgroup. In particular all 
the $2$-torsion of $F\o(B)$ is in $T$. In $F\o(B)$ modulo $T$ one also has 
$f'^{2}=f^{2}$ in this quotient. Now in groups of finite Morley rank without 
involutions any element has a unique square root 
(\cite[p. 72]{BorovikNesin(Book)94}, \cite[Fact 2.25]{AltBorCher97}), 
and one concludes that $f=f'$ modulo $T$, i.e., $f'=ft$ for some $t$ in $T$. 
As $T$ is central in $F\o(B)$, one gets by taking squares that $t^{2}=1$, 
and as $T$ has only finitely many involutions it implies that there are only finitely many 
possibilities for $f'$, once $f$ is fixed. 
This gives the desired finiteness of the fibers of the multiplication map $P$. 

So our statement reduces to proving that the map $P$ is onto. 

In connected solvable groups $B$ of finite Morley rank, Carter subgroups 
are conjugate and cover the connected nilpotent quotient $B/F\o(B)$ 
by \cite{Nesin90-b} and 
\cite[Theorem 3.11 and Proposition 5.1]{FreconJaligot07}. Combined with 
Lemma \ref{LempToralCarter} this yields 
$$B=C\o(j)F\o(B)$$ 
and for our desired factorisation it suffices thus to show that 
$$F\o(B)=C\o_{F\o(B)}(j)\cdot [F\o(B)]^{-_{j}}.\leqno{(\dagger)}$$ 
But the Sylow $2$-subgroup of $F\o(B)$ is a central $2$-torus, again by 
earlier results of Nesin according to which divisible torsion is central in connected 
nilpotent groups of finite Morley rank. Hence Theorem \ref{actionalpha} can 
be applied with the action by conjugation of $j$ on $F\o(B)$, and this 
gives exactly the desired factorisation $(\dagger)$ as above. 

$(2)$. 
If $c \in C\o(j)$ and $u\in U^{-_{j}}$, then 
$(u^{c})^{j}=u^{cj}=u^{jc}=u^{-c}=(u^{c})^{-1}$, and as $u^{c}\in U$ by normality 
of $U$ in $B$ one gets $u^{c}\in U^{-_{j}}$. Hence the subgroup 
$U^{-_{j}}$ is normalized by $C\o(j)$, and as it is also central in $F\o(B)$ it is 
normal in $B=C\o(j)F\o(B)$. 
\qed

\bigskip
The following application of Lemma \ref{computationinside} will be 
used in the most critical situations in \cite{DeloroJaligotII}. 
In the second claim of the corollary below the two involutions considered are not 
necessarily conjugate, a point we will use fully in the analysis 
of \cite{DeloroJaligotII}. 

\begin{cor}
Let $G$ be a group of finite Morley rank, $i$ an involution of $G$, 
and assume $C\o(i)\leq B$ for some definable connected 
solvable subgroup $B$ of odd type. Let $j$ be an involution of $B$. 
Then 
\begin{itemize}
\item[$(1)$]
$\rk(B)=\rk(C\o_{B}(j))+\rk([F\o(B)]^{-_j})$. 
\item[$(2)$]
In particular, if $C\o(j)\leq B$ and $C\o(i)$ and $C\o(j)$ have the same rank 
(which is the case in particular if $i$ and $j$ are conjugate), then 
$\rk ([F\o(B)]^{-_j})=\rk([F\o(B)]^{-_i})$. 
\end{itemize}
\end{cor}
\proof
By assumption, $C\o(i)$ contains trivial $2$-unipotent subgroups. It follows 
from the torality result of \cite{BurdgesCherlinSemisimpleTorsion} that $i$ belongs 
to a $2$-torus, say $T$. Then $i\in T\leq C\o(i)\leq B$. Now 
Lemma \ref{computationinside} applies in $B$. 
\qed

We also mention the following additional information to Lemma \ref{computationinside} 
when the Pr\"{u}fer $2$-rank is $1$. Of course, this corresponds more and 
more to abstract descriptions of the Borel subgroups of $\PSL_{2}$ or $\SL_{2}$ 
over algebraically closed fields of characteristic different from $2$. 

\begin{lem}\label{CompInsidePruef1}
Adopt the same assumptions and notations as in Lemma \ref{computationinside} 
and assume furthermore that the Pr\"{u}fer $2$-rank is $1$. 
Then 
\begin{itemize}
\item[$(1)$]
$C\o(j)<B$ if and only if $F(B)$ contains no involutions. 
\item[$(2)$]
$C_{B}(j)$ is connected, $B=C\o(j)\cdot [F\o(B)]^{-_{j}}$, and the multiplication map 
giving this decomposition is one-to-one whenever $C\o(j)<B$. 
\item[$(3)$]
The set of involutions of $B$ is exactly $j^{F\o(B)}$. 
\end{itemize}
\end{lem}
\proof
First, we note that a connected solvable group of finite Morley rank has connected 
Sylow $2$-subgroups. In particular a connected solvable group of finite Morley rank of 
odd type and of Pr\"{u}fer $2$-rank $1$ has Sylow $2$-subgroups isomorphic to 
$\Z_{2^{\infty}}$. 

$(1)$.
If $F\o(B)$ contains an involution, its maximal $2$-torus is nontrivial. 
As it is central in $B$ by 
\cite[Fact 2.12 (1)]{DeloroJaligotII}, 
the toral involution is central in $B$, and it must be $j$. 
If the Sylow $2$-subgroup of $F(B)$ is finite, then it is central in $B$ by 
\cite[Fact 1.2]{DeloroJaligotII}. 
One gets again that an involution of $B$ must be central in $B$ and by 
connectedness of the Sylow $2$-subgroups of $B$ it must be $j$ again. 

Conversely, if $C_{B}(j)=B$, then $F(B)$ contains the involution $j$. 

$(2)$. 
We can prove the factorisation much more directly here. By $(1)$, we may assume 
$F\o(B)$ without involutions. 
We work in $B/F\o(B)$. As this quotient is abelian by 
\cite{Nesin90-b}, 
$j$ induces by conjugation a trivial action on this quotient. 

In particular, for every $b\in B$, there exists $f\in F\o(B)$ such that 
$b^{j}=bf$. Conjugating again by $j$ one gets $f^{j}=f^{-1}$, so 
$f\in [F\o(B)]^{-_{j}}$. By $(1)$ and Fact \ref{ActionInvOn2PerpGroup}, 
$f=g^{2}$ for some $g$ still in $[F\o(B)]^{-_{j}}$. Then 
$b=(bg)g^{-1}$, where $(bg)^{j}=(bf)g^{-1}=bg\in C_{B}(j)$ and 
$g\in [F\o(B)]^{-_{j}}$.

Hence $B=C_{B}(j)\cdot [F\o(B)]^{-_{j}}$. 
We show the uniqueness of this decomposition whener $C\o(j)<B$. 
If $c_{1}f_{1}=c_{2}f_{2}$, with natural notations, then $f_{2}f_{1}^{-1}\in C(j)$, and 
$f_{1}^{2}=f_{2}^{2}$, which in the group without involutions 
$F\o(B)$ implies $f_{1}=f_{2}$, and then the uniqueness of the decomposition follows. 

The product map corresponding to the decomposition 
$B=C_{B}(j)\cdot [F\o(B)]^{-_{j}}$ 
has finite fibers. Hence each definable generic subset of the source set has a 
generic image in $B$. When the product map is one-to-one, this forces 
that $C_{B}(j)$ has Morley degree one. Otherwise, $C_{B}(j)=B$ is also connected. 

$(3)$. 
As the Pr\"{u}fer $2$-rank is $1$, all involutions of the connected solvable group $B$ 
are conjugate by the structure of Sylow $2$-subgroups. 
Then $(2)$ gives the desired equality. 
\qed

\subsection{Groups of Pr\"{u}fer $2$-rank $1$}\label{SectionGpsPruferRankOne}

It is proved in \cite[Lemma 2.34]{CherlinJaligot2004} that if $S$ is a 
Sylow $2$-subgroup of a group of finite Morley rank of odd type and of 
Pr\"{u}fer $2$-rank $1$, and if $w$ is an involution conjugate to the unique 
involution $i$ of $S\o$, then $w$ cannot centralize $S\o$ unless it is $i$. From this 
one can deduce that such involutions $w$ which normalize $S\o$ but are not 
inside must invert $S\o$, and one gets in this situation a subgroup of the Sylow 
$2$-subgroup isomorphic to that of Chevalley groups of type $\PSL_{2}$ over 
algebraically closed fields of characteristic different from $2$. 

We are going to delineate entirely the structure of Sylow $2$-subgroups of connected 
groups of odd type and Pr\"{u}fer $2$-rank $1$, showing in any case analogy 
with the structure of Sylow $2$-subgroups in groups 
of the form $\PSL_{2}$ or $\SL_{2}$ over some algebraically closed 
field of characteristic different from $2$ or of connected solvable groups. 
First, we recall from the main result of \cite{BurdgesCherlinSemisimpleTorsion} 
that any $2$-element of a connected group of finite Morley rank of odd type 
is {\em $2$-toral}, i.e., belongs to a $2$-torus. When 
the Pr\"{u}fer $2$-rank is $1$ it implies by Fact \ref{FactStruct2Syl} that all 
involutions are conjugate. 

\begin{prop}\label{PropStructSPruf1}
Let $G$ be a connected group of finite Morley rank of odd type and of Pr\"{u}fer 
$2$-rank $1$. Then there are exactly three possibilities for the isomorphism type of 
a Sylow $2$-subgroup $S$ of $G$. 
\begin{enumerate}
\item[$(1)$]
$S=S\o\rtimes{\< w \>}$ for some involution 
$w$ which acts on $S\o$ by inversion. 
\item[$(2)$]
$S=S\o\cdot{\< w \>}$ for some element 
$w$ of order $4$ which acts on $S\o$ by inversion.
\item[$(3)$]
$S=S\o$.
\end{enumerate}
\end{prop}

Our proof of Proposition \ref{PropStructSPruf1} uses the following results. 

\begin{fait}\label{C(So)=So}
{\bf \cite{BorovikCherlin08}}
Let $G$ be a connected group of finite Morley rank of odd type, and fix 
some Sylow $2$-subgroup $S$ of $G$. Then $C_{S}(S\o) = S\o$.
\end{fait}

\begin{fait}\label{AutZtwoinf}
$\Aut(\Z_{2^{\infty}}) \simeq \Z_{2}^\times$. 
In particular the only nontrivial automorphism 
of finite order of $\Z_{2^{\infty}}$ is the inversion. 
\end{fait}
\proof
It is clear that $\Aut(\Z_{2^{\infty}})$ is isomorphic to the group of invertible 
elements of the ring $\Z_{2}$. But 
${\Z_{2}}^{\times} \simeq {\Z/2\Z \times \Z_{2}}$ and the right factor is torsion-free 
as the characteristic is $0$. In particular $\Aut(\Z_{2^{\infty}})$ has only one 
nontrivial automorphism of finite order, inversion. 
\qed

\begin{lem}\label{lemPruef1}
Let $G$ be a connected group of finite Morley rank of odd type and of 
Pr\"{u}fer $2$-rank $1$, and $S$ a Sylow $2$-subgroup of $G$. 
Then $[S:S\o] \leq 2$ and elements of $S\setminus S\o$ act on $S\o$ by inversion.
\end{lem}
\proof
By Fact \ref{C(So)=So}, every element of $S\setminus S\o$ has a nontrivial 
action on $S\o$. (This is the original argument in \cite[Lemma 2.34]{CherlinJaligot2004} 
in the light of the torality of \cite{BurdgesCherlinSemisimpleTorsion}.) 
 
Now Fact \ref{AutZtwoinf} implies that this action is by inversion, 
and in particular all elements of $S\setminus S\o$ have the same action on $S\o$. 
By Fact \ref{C(So)=So} again one gets that there is at most one coset of 
$S\o$ distinct from $S\o$ in $S$, which proves our claim. 
\qed

\bigskip
\noindent
{\bf Proof of Proposition \ref{PropStructSPruf1}.} 
Assume that $S$ is not connected. We prove that $S$ is either 
isomorphic to the Sylow $2$-subgroup of $\PSL_2$ or to that of $\SL_2$ in characteristic 
different from $2$.

Let $w \in S\setminus S\o$. By Lemma \ref{lemPruef1}, 
$S = S\o \cdot \< w \>$, $w$ inverts $S\o$, and $w^2 \in S\o$.
If $w^2 = 1$, then $S$ is obviously as in $\PSL_2$. 
Now suppose $w^2\neq 1$. As $w$ inverts $S\o$ and 
$w^2 \in S\o$, one has $w^2 = (w^2)^w = (w^2)^{-1}$. 
So $w^2$ is the involution of $S\o$, as in $\SL_2$.
\qed

\bigskip
Remarkably, the conclusion of Proposition \ref{PropStructSPruf1} will be obtained 
in a very general setting in the case of connected locally$\o$ solvable$\o$ groups of 
odd type, in the early analysis of \cite{DeloroJaligotII}. 
This has the effect of simplifying slightly 
certain arguments in the entire classification of small groups of odd type, not 
spectacularly as the main difficulties with strongly real elements remain throughout, 
but at least morally and technically, removing certain residual involutions sometimes 
occuring in Weyl groups. As an example of such minor but seemingly 
simplifications, it entirely eliminates the need for the lengthy 
Sections 6.1 and 6.2 in the analysis of Weyl groups in \cite{CherlinJaligot2004}. 

\subsection{The Borovik-Cartan decomposition}\label{SectionBorovikCartan}

In quadratic algebra, the {\em Cartan polar form} 
yields a decomposition $f=u\cdot s$ of any automorphism $f$ of $\C^{n}$,
with $u$ in the unitary group and $s$ a self-adjoint automorphism. 
It is well known by french {\em taupins} that the Cartan decomposition can be
used to prove connectedness of the unitary group. Lines of arguments, formally 
similar, were implicitely (and unconsciously) rediscovered to prove connectedness 
of centralizers of involutions in 
certain context of groups of finite Morley rank in \cite{BorovikBurdgesCherlin07}. 
We are going to rework this in a form more suitable for 
\cite{DeloroJaligotII}, but let us first serve some refreshments. 

We denote (Hermitian) adjunction on $\C^{n}$ by ${}^{*}$, so that 
$$(xy)^{*}=y^{*}x^{*} \mbox{~and~} x^{**}=x$$ 
for any automorphisms $x$ and $y$ of $\C^{n}$. 
We note that it follows from the first equality that $1^{*}=1$ and 
$(x^{-1})^{*}=(x^{*})^{-1}$ for any automorphism $x$. 

If $g$ is an automorphism of $\C^{n}$, then $g^{*}g$ is a self-adjoint automorphism, 
and one can sometimes get a self-adjoint square root $s$ of $g^{*}g$.
Letting $u=gs^{-1}$, one finds 
$$u^{*}u = (gs^{-1})^{*}(gs^{-1}) = s^{-1}g^{*}gs^{-1} = 1$$ 
as $g^{*}g=s^{2}$, and thus $u$ is orthogonal with respect to $*$ 
(which might be called unitary here). 
Moreover, it is usually possible to choose $s$ in a consistent fashion when $g$ varies, 
and the decomposition turns out to be continuous.
Using the assignment $g \mapsto u$ to map $\GL_{n}$ onto $U_{n}$, one finds 
the connectedness of the latter.

An analog of the Cartan decomposition was used implicitly in \cite{BorovikBurdgesCherlin07} with 
$${g}^{*} = g^{-i}$$ 
for arbitrary elements $g$ and a fixed involution $i$ of a given group $G$ of finite Morley 
rank. This operation on elements of the (ambient) group behaves formally like 
an adjunction as one sees immediately that $(gg')^{*}=g'^{*}g^{*}$ 
for any elements $g$ and 
$g'$ and that $1^{**}=1$. Using this notation, one notices also that 
$$g^{*}g=ig^{-1}ig=ii^{g},$$
i.e., $g^{*}g$ is nothing else than a strongly real element $ii^{g}$, 
inverted by the involution $i$ and its conjugate $i^{g}$. 
We further remark that with our definition the unitary group ``$U(*)$" 
corresponds naturally to $C(i)$ and that the set of self-adjoint automorphisms 
corresponds to the set of strongly real elements inverted by $i$. 
In our abstract context, the only technical point consists thus 
in finding a well defined, and definable, square root function corresponding to the 
extraction of square roots of $ii^{g}=g^{*}g$. 

We prove the following theorem using this decomposition. 

\begin{theo}\label{LemCommutInvolEquivCiCon}
Let $G$ be a connected group of finite Morley rank in which commuting is an 
equivalence relation on the set of involutions, and with no nontrivial normal 
$2$-unipotent subgroup. Then $C(i)$ is connected for any involution $i$ of $G$. 
\end{theo}

Before passing to the proof of Theorem \ref{LemCommutInvolEquivCiCon}, it is worth 
commenting on the assumption of the absence of a normal $2$-unipotent subgroup. 
In fact, the conclusion of Theorem \ref{LemCommutInvolEquivCiCon} 
may fail if one drops this assumption, as the 
following example shows. If $K$ is an algebraically closed field of characteristic $2$, 
then in the connected solvable group of matrices 
$${\left\{
 {\left(
      \begin{array}{llll} t & a \\
      0 & t^{4} \end{array}
 \right)}
\mbox{ : } t \in K^{\times} \mbox{ , } a \in K_{+} 
\right\}}$$
one sees that the centralizer of every involution is the cyclic extension of order $3$ of the 
$2$-unipotent subgroup (consisting of strictly upper triangular matrices). In particular 
centralizers of involutions are not connected in this group, even though commutation 
is an equivalence relation on the set of involutions. 

\begin{lem}\label{LemmaProd2UnipNormal} 
Let $G$ be a connected group of finite Morley rank in which commutation is an 
equivalence relation on the set of involutions, and assume there exists a 
noncentral involution $i$ of $G$. Then, either 
\begin{itemize}
\item[$(1)$]
The definable subset $X$ of elements $g$ of $G$ such that $Z(C(ii^{g}))$ contains no 
involutions is generic in $G$, or 
\item[$(2)$]
$i^{G}$ generates an infinite normal elementary abelian $2$-group. 
\end{itemize}
\end{lem}
\proof
Assuming we are not in case $(1)$, then the 
definable subset $Y$ of elements $g$ such that $Z(C(ii^{g}))$ contains an involution, that is 
the complement of $X$, is generic in $G$. Notice 
that $Y$ is a union of right cosets of $C(i)$. 

For any $g$ in $Y$, $Z(C(ii^{g}))$ 
contains a nontrivial Sylow $2$-subgroup, normalized by $i$ and $i^{g}$. Hence 
each involution $i$ and $i^{g}$ has to commute with an involution in 
$Z(C(ii^{g}))$, and the transitivity of the commutation of involutions implies that 
$i$ and $i^{g}$ commute. Now the transitivity of the commutation of involutions 
forces that any two involutions $i^{g_{1}}$ and $i^{g_{2}}$, where $g_{1}$ and 
$g_{2}\in Y$, have to commute also. In particular $i^{Y}$ is a set of 
commuting involutions. 
As $Y$ is generic in $G$, $i^{Y}$ is generic in $i^{G}$. 

Now $i^{G}$ is infinite as $C(i)<G$ by assumption, 
and $i^{G}$ has Morley degree one as it is in definable bijection 
with $G/C(i)$. As $i^{Y}$ is generic in $i^{G}$, it follows that any two 
$G$-conjugates of $i^{Y}$ intersect on a subset generic in $i^{G}$, and thus 
in particular any two such conjugates, consisting of pairwise commuting involutions, 
commute by the transitivity of the commutation on the set of involutions. 
As $G$ acts transitively by conjugation on $i^{G}$, it follows that $i^{G}$ consists of 
pairwise commuting involutions. But $i^{G}$ is a normal subset of $G$. Hence we conclude 
that $i^{G}$ generates an infinite elementary abelian $2$-group, normal in $G$. 
We are thus necessarily in case $(2)$. 

It remains just to show that cases $(1)$ and $(2)$ are mutually exclusive. But in 
case $(2)$, one has that $ii^{g}$ is an involution, necessarily in $Z(C(ii^{g}))$, 
for any $g$ in $G\setminus C(i)$. Hence $X\subseteq C(i)$, and $X$ is in particular 
not generic in $G$ 
\qed

\bigskip
\noindent
{\bf Proof of Theorem \ref{LemCommutInvolEquivCiCon}.}
If $i$ belongs to the center of $G$, then $C(i)=G$ is connected. 
Assume now $C(i)<G$. By Lemma \ref{LemmaProd2UnipNormal}, 
the definable subset $X$ of elements $g$ of $G$ such that $Z(C(ii^{g}))$ 
contains no involutions is generic. Notice that $X$ is a union of right cosets of $C(i)$. 
We can define the Borovik-Cartan function associated to the operation 
$$g^{*}=g^{-i}.$$

For $g$ in $X$, there exists a unique square root $s$ of $ii^{g}=g^{*}g$ in 
$Z(C(ii^{g}))$, and the function 
$$\psi~:~g\mapsto gs^{-1}$$ 
from $X$ to $G$ is well defined and definable. We also note that the square root $s$ is 
necessarily in the subgroup of $Z(C(ii^{g}))$ of elements inverted by $i$, as this 
subgroup has no involutions also, and thus $s^{*}=s$ and $(s^{-1})^{*}=s^{-1}$. 
We then get 
$$\psi(g)^{*}\psi(g)=(gs^{-1})^{*}gs^{-1}=s^{-1}g^{*}gs^{-1}=1$$ 
as $g^{*}g=s^{2}$. In other words, the function $\psi$ takes all its values in $C(i)$. 

Furthermore, for any $c$ in $C(i)$ one has 
$$(cg)^{*}(cg)=g^{*}c^{*}cg=g^{*}g,$$ 
so that the square roots $s_{cg}$ and $s_{g}$ of $(cg)^{*}$ and $g^{*}$ are the same. 
This shows the following covariant property: 
$$\psi(cg)=(cg)s_{cg}^{-1}=(cg)s_{g}^{-1}=c(gs_{g}^{-1})=c\psi(g).$$
In particular, the fibers of $\psi$ are of constant rank, say $f$, and any subset of $C(i)$ of rank $r$ 
lifts to a subset of $X$ of rank $r+f$. If $C(i)$ were not connected, then it would have two 
disjoint generic subsets of full rank, which would necessarily lift to disjoint generic 
subsets of $X$, contradicting the connectedness of $G$. Hence $C(i)$ is connected. 
\qed

\bigskip
In connection of Theorem \ref{LemCommutInvolEquivCiCon}, it is worth mentioning 
the so-called {\em $Z^{*}$-theorem} (correcting also a slightly inacurate statement in 
\cite{BorovikBurdgesCherlin07}). We say that an involution $k$ of a group $G$ of 
finite Morley rank is {\em isolated} whenever 
$$|k^{G}\cap S|=1$$ 
for some (any) Sylow $2$-subgroup $S$ of $G$. 

\begin{fait}\label{Zstar}
{\bf ($Z^{*}$-theorem \cite[Theorem 6]{BorovikBurdgesCherlin07})}
Let $G$ be a connected group of finite Morley rank. 
If some involution $k$ of $G$ is isolated, then $C(k)$ is connected.
\end{fait}

We stress the fact that one only has an implication as in Fact \ref{Zstar}, 
and not an alternative ``either $k$ is not isolated or $C(k)$ is connected" as suggested 
by the statement used in \cite[Theorem 6]{BorovikBurdgesCherlin07}. In fact there 
are connected groups of finite Morley rank with non-isolated involutions $k$ and with 
$C(k)$ connected, as the following examples show. 
In $\SL_{2}$ over some algebraically closed field of characteristic $2$, 
centralizers of involutions are connected and all involutions of the elementary 
abelian $2$-subgroups are conjugate at the same time. Another example 
is provided by $\GL_{2}$ over an algebraically closed field of characteristic different from 
$2$, where the centralizer of the involution 
$$
{\left(
      \begin{array}{llll} 1 & 0 \\
      0 & -1 \end{array}
 \right)}
$$
is the definable subgroup of diagonal matrices (a connected torus), and this involution 
is conjugate to the other non-central involution 
$$
{\left(
      \begin{array}{llll} -1 & 0 \\
      0 & 1 \end{array}
 \right)}.
$$

To conclude, we mention the following special case of 
Theorem \ref{LemCommutInvolEquivCiCon}, which will be used notably 
in our reworking of \cite[Case I]{BurdgesCherlinJaligot07} in 
\cite{DeloroJaligotII}. 

\begin{cor}\label{CorFinal}
Let $G$ be a connected group of finite Morley rank such that 
$U_{2}(F(G))=1$. 
\begin{itemize}
\item[$(1)$]
If $G$ contains a strongly embedded subgroup in which an involution is central, 
then the centralizers of involutions of $G$ are connected. 
\item[$(2)$]
In particular, if $G$ has odd type, Pr\"{u}fer $2$-rank $1$, and connected 
Sylow $2$-subgroups, then the centralizers of involutions of $G$ are connected. 
\end{itemize}
\end{cor}
\proof
$(1)$. Let $M$ denote the strongly embedded subgroup of $G$. By assumption, 
$M=C(i)$ for some of its involutions, and as all involution of $M$ are 
$M$-conjugate (\cite[Theorem 10.19]{BorovikNesin(Book)94}) 
this holds for any involution of $M$. We also note that the set of involutions 
of $M$ form an elementary abelian $2$-group (central in $M$). 
By strong embedding of $M$, we also see that any two distinct conjugates of the elementary 
abelian $2$-subgroup of $M$ have a trivial intersection. Hence 
Theorem \ref{LemCommutInvolEquivCiCon} gives the desired result. 

$(2)$. By assumption, a Sylow $2$-subgroup $S$ of $G$ has the isomorphism type 
$$S\simeq \Z_{2^{\infty}}.$$ 
Let $i$ denote the unique involution of $G$, and $M=C(i)$. As $i$ is the unique 
involution of $M$, $M$ is strongly embedded in $G$. Hence we are in a special case 
of case $(1)$. 
\qed

\bigskip
We note that Corollary \ref{CorFinal} $(2)$ follows also more generally 
from the recent proof that $C(i)/C\o(i)$ has exponent at most $2$ for any involution 
$i$ in a connected group of finite Morley rank of odd type \cite{DeloroBagatelle}. 
Actually, this result implies, more generally, that centralizers of involutions are connected 
in any connected groups of finite Morley rank of odd type with connected 
Sylow $2$-subgroups.

\bibliographystyle{alpha}
\bibliography{biblio}

\end{document}